\documentclass{amsart}
\numberwithin{equation}{section}
\usepackage{amscd}
\usepackage{amssymb}
\let\cal\mathcal

\def\Oscr{{\cal O}}

\let\blb\mathbb

\def\GG{{\blb G}}

\let\at\ast
\def\Der{\operatorname{Der}}

\def\Lie{\mathop{\text{Lie}}}

\def\Supp{\mathop{\text{\upshape Supp}}}

\def\Spec{\operatorname {Spec}}

\def\GL{\operatorname {GL}}

\def\End{\operatorname {End}}

\def\im{\operatorname {im}}

\def\coker{\operatorname {coker}}
\def\ker{\operatorname {ker}}

\def\End{\operatorname {End}}

\def\rk{\operatorname {rk}}

\def\r{\rightarrow}

\newtheorem{lemma}{Lemma}[section]
\newtheorem{proposition}[lemma]{Proposition}
\newtheorem{theorem}[lemma]{Theorem}
\newtheorem{corollary}[lemma]{Corollary}

\theoremstyle{definition}

\newtheorem{example}[lemma]{Example}

{

}

\theoremstyle{remark}

\newtheorem{remark}[lemma]{Remark}

\newdimen\uboxsep \uboxsep=1ex
\def\uboxn#1{\vtop to 0pt{\hrule height 0pt depth 0pt\vskip\uboxsep
\hbox to 0pt{\hss #1\hss}\vss}}

\def\uboxs#1{\vbox to 0pt{\vss\hbox to 0pt{\hss #1\hss}
\vskip\uboxsep\hrule height 0pt depth 0pt}}

\def\codim{\operatorname{codim}}
\def\si{\operatorname{si}}
\def\trdeg{\operatorname{trdeg}}

\def\ad{\operatorname{ad}}
\def\sm{\operatorname{sm}}
\title[A degree inequality]{A degree inequality for Lie algebras with 
a regular Poisson semi-center}
\author{A. I. Ooms}
\author{M. Van den Bergh}
  \email{alfons.ooms@uhasselt.be, michel.vandenbergh@uhasselt.be}
\address{Mathematics Department, Hasselt University, Agoralaan, 3590 Diepenbeek, Belgium.}
\thanks{The second author is a senior researcher at the FWO}
\begin{document}
\begin{abstract} For Lie algebras whose Poisson semi-center 
is a polynomial ring we give a bound for the sum of the degrees of the
  generating semi-invariants.  This bound was previously known in many
special cases.
\end{abstract}\maketitle
\section{Introduction}
In this paper we work over an algebraically closed base field $k$ of
characteristic zero.  Let $\frak{g}$ be a finite dimensional Lie
algebra.  A non-zero element $f\in S\frak{g}$ is called a
\emph{semi-invariant with weight $\chi \in \frak{g}^\ast$} if for all
$v\in \frak{g}$ we have
\[
\ad(v)(f)=\chi(v)f
\]
We say that a semi-invariant is \emph{proper} if $\chi\neq 0$.  The
$k$-algebra generated by the semi-invariants in $S\frak{g}$ is denoted
by $(S\frak{g})^{\frak{g}}_{\si}$. This ring is called the
\emph{Poisson semi-center} of $S\frak{g}$.

The \emph{stabilizer} of $x\in \frak{g}^\ast$ is denoted by $\frak{g}^\ast_x$. I.e.
$ \frak{g}_x=\{v\in \frak{g}\mid \forall w\in \frak{g}:x([v,w])=0\} $.
The minimal value of $\dim \frak{g}_x$ is called the \emph{index} of
$\frak{g}$ and is denoted by~$i(\frak{g})$. 
An element $x\in
\frak{g}^\ast$ is called \emph{regular} if $\dim
\frak{g}_x^\ast=i(\frak{g})$. The regular elements form an open dense
subset of $\frak{g}^\ast$ which we denote by $\frak{g}_{\text{reg}}^\ast$.

The following is our main result. 
\begin{theorem} \label{ref-1.1-0} (see Prop.\ \ref{ref-3.1-10} and Prop.\
  \ref{ref-5.7-39} below.) Assume that $(S\frak{g})^{\frak{g}}_{\si}$ is
  freely generated by homogeneous elements
  $f_1,\ldots,f_r$. Then
\begin{equation}
\label{ref-1.1-1}
\sum_{i=1}^r\deg f_i\le \frac{1}{2}\bigl(\dim \frak{g}+i(\frak{g})\bigr)
\end{equation}
\end{theorem}
It is well-known that \eqref{ref-1.1-1} holds for semi-simple Lie
algebras \cite[Thm.\ 7.3.8]{Dixmier} and Frobenius Lie algebras
\cite[pp.\ 339--343]{Ooms3}.  Numerous other special cases are known
(e.g.\ \cite{JFM,Ooms2}). 

\medskip

For Theorem \ref{ref-1.1-0} to be valid in the stated generality it is
essential that we consider semi-invariants instead of invariants, as
the following trivial example by Panyushev shows.
\begin{example} Let $\frak{g}=kv_1+kv_2+kv_3+kv_4$ with non-trivial
  brackets $[v_1,v_2]=v_2$, $[v_1,v_3]=v_3$, $[v_1,v_4]=-v_4$. Then
  $\dim \frak{g}=4$, $i(\frak{g})=2$. The generating invariants are
  $v_2v_4$ and $v_3v_4$. So the sum of their degrees is $2+2=4$ which
  is strictly bigger than $1/2(\dim \frak{g}+i(\frak{g}))=3$. However
  the generating \emph{semi-invariants} are $v_2,v_3,v_4$ and the
  sum of their degrees is 3, which does not violate the inequality.
\end{example}
For brevity we will call a Lie algebra
\emph{coregular} if $(S\frak{g})^{\frak{g}}$ is a polynomial ring. 
\begin{corollary} Assume that 
  $(S\frak{g})^{\frak{g}}_{\si}=(S\frak{g})^{\frak{g}}$ and $\frak{g}$ is
coregular with center $Z(\frak{g})$.  Then
\begin{equation}
\label{ref-1.2-2}
3\,i(\frak{g})\le \dim \frak{g}+2\dim Z(\frak{g})
\end{equation}
\end{corollary}
\begin{proof}
In this situation we have the equality $r=i(\frak{g})$ 
(see Proposition \ref{ref-4.1-16} below). The observation
that $\deg f_i\ge 2$, unless $f_i\in
Z(\frak{g})$ yields
\[
\dim Z(\frak{g})+2(i(\frak{g})-\dim Z(\frak{g}))\le \sum_{i=1}^r \deg f_i
\le \frac{1}{2}(\dim \frak{g}+i(\frak{g}))
\]
which translates into \eqref{ref-1.2-2}.
\end{proof}
 The number on the right hand side of \eqref{ref-1.1-1} occurs frequently
 in the theory of enveloping algebras. 
For example it is an upper bound for the transcendence degree of a maximal 
commutative subfield of the division ring of fractions of $U\frak{g}$
and this bound can be achieved in many cases \cite{Joseph:gk, Ooms5}. 
Likewise by a result of Sadetov \cite{Sadetov} it is
the  maximum 
 transcendence degree of a Poisson commutative subfield of the field
 of fractions of~$S\frak{g}$.

\medskip

For the proof of Theorem \ref{ref-1.1-0} we first reduce
to the case that there are no proper semi-invariants (i.e.\
  $(S\frak{g})^{\frak{g}}_{\si}=(S\frak{g})^{\frak{g}}$).  
  In this situation one may prove a result which is more precise than
  Theorem~\ref{ref-1.1-0}.  Assume first that $\frak{g}$ is
  non-abelian. Let $B=([v_i,v_j])_{ij}\in M_{n}(S\frak{g})$ be the
  \emph{structure matrix} of $\frak{g}$ where $v_1,\ldots, v_n$ is an
  arbitrary basis of~$\frak{g}$. Put $s=\dim \frak{g}-i(\frak{g})$.
  Then the greatest common divisor of the $s\times s$-minors in $B$ is
  a semi-invariant in $S\frak{g}$~\cite{Ooms3}. Below we will call it
  the \emph{fundamental semi-invariant} and we denote its degree by
  $d(\frak{g})$.  If $\frak{g}$ is abelian we put $d(\frak{g})=0$.

\begin{proposition} (see Prop.\ \ref{ref-5.7-39}.) \label{ref-1.4-3} Assume that
  $(S\frak{g})^{\frak{g}}_{\si}=(S\frak{g})^{\frak{g}}$ and $\frak{g}$
  is coregular.  Then we have
\begin{equation}
\label{ref-1.3-4}
\sum_{i=1}^r\deg f_i=\frac{1}{2}\bigl(\dim
\frak{g}+i(\frak{g})-d(\frak{g})\bigr)
\end{equation}
\end{proposition}
Taking into account Propositions \ref{ref-4.1-16} and \ref{ref-5.1-26}
below, this result may also be deduced from \cite[Remark
1.6.3]{Panyushev}.  Our proof uses the general techniques from
\cite{Knop3,Knop1} and is quite different from \cite{Panyushev}. We
obtain a certain nice complex of length three, consisting of free
$S\frak{g}$-modules which, besides implying \eqref{ref-1.3-4}, yields
some additional information on $\frak{g}^\ast\setminus \frak{g}^{\ast}_{\text{reg}}$ (see Proposition \ref{ref-1.6-6} below).
\begin{corollary}
\label{ref-1.5-5}
 Assume that
  $(S\frak{g})^{\frak{g}}_{\si}=(S\frak{g})^{\frak{g}}$ and $\frak{g}$
  is coregular. Then \eqref{ref-1.1-1} is an equality if and only
if $\frak{g}^\ast\setminus \frak{g}^\ast_{\text{reg}}$ has codimension $\ge 2$.
\end{corollary}
This follows from the easily verified fact that $d(\frak{g})=0$ if and
only if $\codim_{\frak{g}^\ast}(\frak{g}^\ast\setminus
\frak{g}^\ast_{\text{reg}}) \ge 2$.

Proposition \ref{ref-1.4-3} is false without the assumption
$(S\frak{g})^{\frak{g}}_{\si}=(S\frak{g})^{\frak{g}}$. Counter
examples are given by Frobenius Lie algebras. By definition these satisfy
$i(\frak{g})=0$ and thus the fundamental semi-invariant is equal to
$\det B$. Hence $d(\frak{g})=\dim \frak{g}$ and the righthand side of
\eqref{ref-1.3-4} is zero. Since $(S\frak{g})^{\frak{g}}_{\si}$ is
freely generated by the irreducible factors of $\det B$ \cite{Ooms3}
the lefthand side of \eqref{ref-1.3-4} is never zero. It would be
interesting to find a version of Proposition \ref{ref-1.4-3} which
holds in the same generality as Theorem~\ref{ref-1.1-0}.

\medskip

As mentioned above we may use our methods to obtain some additional
necessary conditions for coregularity.  As $\frak{g}$ acts by
derivations on $S\frak{g}$ we have a $S\frak{g}$-linear map
\[
\rho:S\frak{g}\otimes \frak{g}\r
\Der_k(S\frak{g})=S\frak{g}\otimes \frak{g}^\ast
\]
\begin{proposition} \label{ref-1.6-6} (see Prop.\ \ref{ref-5.4-34},\ref{ref-5.10-42}) Assume that $(S\frak{g})^{\mathfrak{g}}_{\si}=(S\frak{g})^{\mathfrak{g}}$
  and $\frak{g}$ is coregular. Then
\begin{enumerate}
\item
$\ker \rho$ is a free $S\frak{g}$-module;
\item if $\frak{g}$ is not abelian then
$
\codim(\frak{g}^\ast\setminus\frak{g}^\ast_{\text{reg}})
\le 3
$.
\item If $\codim(\frak{g}^\ast\setminus\frak{g}^\ast_{\text{reg}})=3$ then 
$\frak{g}^\ast\setminus\frak{g}^\ast_{\text{reg}}$ is purely of codimension
three. 
\end{enumerate}
\end{proposition}

\begin{example} 
\label{ref-1.7-7} We illustrate the above results with an easy example.
  For $n\ge 3$ let $\frak{g}=L(n)$ be the $n$-dimensional
  \emph{standard filiform Lie algebra}. $L(n)$ has a basis
  $v_1,\ldots,v_n$ and non-trivial brackets $[v_1,v_i]=v_{i+1}$ for
  $i=2,\ldots,n-1$. In this case
\[
(S\frak{g})^{\frak{g}}=k[v_2,\ldots,v_n]^e
\]
where $e$ is the derivation 
\[
e=\sum_{i=2}^{n-1} v_{i+1}\frac{\partial\ }{\partial v_i}
\]
Dixmier verified by direct computation that $L(3)$, $L(4)$ are
coregular but $L(5)$ is not~\cite{Dixmier2}. From the classical
correspondence between $\GG_a$-invariants and $\operatorname{SL}_2$-covariants (e.g.
\cite[\S33]{GY}) one obtains that $L(n)$ is coregular if and only
if $n< 5$ (see e.g. \cite{Tyc}).

\medskip

In order to apply the criteria given above it is advantageous to
use the structure matrix $B$ which was already introduced.  It is easy to
see that the $S\frak{g}$-linear map $\rho$ is represented by the
matrix~$B$. Furthermore $\frak{g}_x=\ker x(B)$. If we write $
r(\frak{g})$ for the rank of $B$ over the quotient field of
$S\frak{g}$ then $ i(\frak{g})=\dim\frak{g}-r(\frak{g}) $.

\medskip

In the case of $L(n)$ the structure matrix
looks like
\[
\begin{pmatrix}
0      &  v_3   &\cdots & v_n    & 0\\
-v_3   &  0     &\cdots & 0      & 0\\
\vdots & \vdots &       & \vdots & \vdots\\
-v_n   &  0     &\cdots & 0      & 0\\
0      &  0     &\cdots & 0      & 0
\end{pmatrix}
\]
We deduce $i(\frak{g})=n-2$. Furthermore the fundamental
semi-invariant is $1$ unless $n=3$ in which case it is $v_3^2$. Thus
\[
d(\frak{g})=
\begin{cases}
0&\text{if $n> 3$}\\
2&\text{if $n=3$}
\end{cases}
\]
Since
$\frak{g}$ is nilpotent there are no proper semi-invariants.  As
$Z(\frak{g})=kv_n$ the numerical criterion \eqref{ref-1.2-2} for
coregularity becomes
\[
3(n-2)\le n+2
\]
which holds iff $n\le 4$. Hence the non-coregularity of $L(n)$ for $n\ge 5$
is detected by~\eqref{ref-1.2-2}.

We have
\[
\frak{g}^\ast\setminus\frak{g}^\ast_{\text{reg}}=\{x\in \frak{g}^\ast\mid x(v_i)=0\text{ for $i=3,\ldots,n$}\}
\]
Thus $\codim
(\frak{g}^\ast\setminus\frak{g}^\ast_{\text{reg}}) =n-2$ and so the fact that
$L(n)$ is not coregular for $n\ge 6$ is detected by the numerical
criterion Prop.\ \ref{ref-1.6-6}(2).

If $n=5$ then 
\[
\ker \rho=\{(A_1,\ldots,A_5)\in k[v_1,\ldots,v_5]\mid A_2v_3+A_3v_4+A_4v_5=0,
A_1v_3=A_1v_4=A_1v_5=0\}
\]
This kernel is minimally generated by $w_1=(0,v_4,-v_3,0,0)$,
$w_2=(0,0,v_5,-v_4,0)$, $w_3=(0,v_5,0,-v_3,0)$ and $w_4=(0,0,0,0,1)$.
These generators are related by $v_5w_1+v_3w_2-v_4w_3=0$ so
they are not free. Thus the non-coregularity of $L(5)$ is detected by
Prop.\ \ref{ref-1.6-6}(1) but not by \ref{ref-1.6-6}(2).

Let us now consider $n=3$. In this case the generating invariant is $v_3$
and the equality \eqref{ref-1.3-4} becomes $1=(1/2)(3+1-2)$.

Assume $n=4$. Now the generating invariants are $v_4$ and
$v_2v_4-(1/2)v_3^2$. Then \eqref{ref-1.3-4} becomes $1+2=(1/2)(4+2-0)$.

\medskip

Although not directly related to the content of this paper let us
remind the reader that much is known classically about the invariant
theory of $\operatorname{SL}_2$. This may be translated back into
results about $L(n)$. 
 For $\frak{g}=L(7)$ one finds that
$(S\frak{g})^{\frak{g}}$ is minimally generated by 23 elements (see
\cite[\S 115]{GY}). On the other hand the transcendence degree of the
fraction field of $(S\frak{g})^{\frak{g}}$ is only $5$.

We refer to \cite{Ooms2} for explicit generators of $(S\frak{g})^{\frak{g}}$ for
many nilpotent Lie algebras of dimension at most 7.
\end{example}
We wish to thank Alexander Elashvili for many stimulating discussions around
this and related problems. 
\section{Preliminaries}
\label{ref-2-8}
 Throughout $\frak{g}$ is a finite dimensional Lie algebra.  
If $V$ is a
finite dimensional representation of $\frak{g}$ then we denote by
$(SV)^{\frak{g}}_{\si}$ the ring of semi-invariants in $SV$.  Note
that if $f$ is a semi-invariant and $g\in SV$ divides $f$ then $g$ is a semi-invariant
as well.  Thus any semi-invariant in $SV$ is a product of
semi-invariants which are irreducible in $SV$.

If $x\in V^\ast$ then $\partial_x$ is the derivation of $SV$ such that
for $v\in V$ we have $\partial_x(v)=x(v)$.

We equip $S\frak{g}$ with the Kostant-Kirillov Poisson bracket of degree $-1$
\[
\{v_1,v_2\}=[v_1,v_2]\qquad (v_1,v_2\in \frak{g})
\]
If $g\in S\frak{g}$ is a semi-invariant with weight $\chi$ then for
all $f\in S\frak{g}$ we have $\{f,g\}=\partial_\chi(f)g$.  From this
we easily deduce the well-known fact that semi-invariants in
$S\frak{g}$ Poisson commute.

\medskip

It will be convenient to
introduce the \emph{Lie-Rinehart algebra} $SV\otimes
\frak{g}$ \cite{Rinehart}. This is a Lie algebra with Lie bracket
\[
[f\otimes v,g\otimes w]=fv(g)\otimes w-gw(f)\otimes v+fg\otimes [v,w]
\]
Sending $f\otimes v$ to $fv(-)$ defines an $S\frak{g}$-linear Lie algebra
homomorphism
\[
\rho:SV\otimes \frak{g}\r \Der_k(SV)
\]
which is called the \emph{anchor map}.  If $(v_i)_i$ is a basis for $\frak{g}$ then
the kernel of the anchor map is given by the sums $\sum_i c_i\otimes v_i$
such that $\sum_i c_iv_i(w)=0$ for all $w\in V$. Note that this kernel is a Lie ideal (as is any kernel of a homomorphism between Lie algebras). 

If we use the identification $\Der_k(SV)= SV\otimes V^\ast$ and we
choose bases $(v_i)_{i=1}^n$, $(w_j)_{j=1}^m$ for $\frak{g}$ and $V$
then the anchor map is represented with respect to these bases by the
\emph{structure matrix} $(v_i(w_j))_{ji}\in M_{m\times n}(SV)$ of $V$.

For convenience we will write $r(V)$ for the rank of the structure
matrix of $V$ over the field of fractions of $SV$. If $\frak{g}$ is in
doubt we write $r_{\frak{g}}(V)$. There is an open subset
$V^\ast_{\text{reg}}$ of $V^\ast$ such that $x\in V^\ast_{\text{reg}}$ iff
$\dim \frak{g}_x=\dim V-r(V)$ where $\frak{g}_x$ denotes the stabilizer
of $x$. I.e. $ \frak{g}_x=\{v\in \frak{g}\mid \forall w\in V:x([v,w])=0\} $.

The \emph{fundamental semi-invariant} in $SV$ is defined as the
greatest common divisor of the $r(V)\times r(V)$ minors in the
structure matrix of $SV$ \cite{Ooms3}, assuming $r(V)>0$. If $C\subset
V^\ast$ is defined by the zeroes of the fundamental semi-invariant
then we have $C\subset V^\ast\setminus V^\ast_{\text{reg}}$ and the
complement of $C\cup V^\ast_{\text{reg}}$ in $V^\ast$ has codimension
$\ge 2$. We write $d(V)$ for the degree of the fundamental
semi-invariant. If $r(V)=0$ (i.e.\ the action of $\frak{g}$ on $V$
is trivial) then we put $d(V)=0$. We record the following
\begin{lemma} 
\label{ref-2.1-9} If $V=\frak{g}$ then the fundamental semi-invariant
in $S\frak{g}$ is the square of the greatest
common divisor of the Pfaffians of the principal $r(\frak{g})\times
r(\frak{g})$ minors in the structure matrix of $\frak{g}$. 
\end{lemma}
\begin{proof} According to \cite{Heymans} any $r(V)\times r(V)$-minor can be
expressed as a quadratic form in Pfaffians of principal $r(V)\times
r(V)$-minors. From this one easily deduces the stated result. 
\end{proof}
In case $V=\frak{g}$ is the adjoint representation
then $r_{\frak{g}}(\frak{g})=r(\frak{g})$ is an even number and we
we have
$r(\frak{g})=\dim
\frak{g}-i(\frak{g})$. 
We put
\begin{align*}
c(\frak{g})&=\frac{1}{2}(\dim \frak{g}+i(\frak{g}))\\
&=\dim\frak{g}-\frac{1}{2}r(\frak{g})
\end{align*}
\section{Reduction to the case without proper semi-invariants}
The following result which generalizes \cite[Thm.\ 1.19(3)]{Ooms3} is
the main result of this section.
\begin{proposition} \label{ref-3.1-10}
Let $\frak{g}$ be a finite dimensional Lie algebra.
Then there exists another finite dimensional Lie
  algebra $\frak{g}'$ such that $(S\frak{g})^{\frak{g}}_{\si}=
(S\frak{g}')^{\frak{g}'}_{\si}=(S\frak{g}')^{\frak{g}'}$.
  Moreover $c(\frak{g}')=c(\frak{g})$.
\end{proposition}
If $\frak{g}$ is almost algebraic then we may take $\frak{g}'$ to be
the intersection of the kernels of the non-trivial weights of the
semi-invariants in $S\frak{g}$ \cite{BGR,Ooms3,JFM,RV}.  This procedure must be
modified for non-almost algebraic Lie algebras.
\begin{example} Let $\frak{g}=kv_1+kv_2+kv_3$ be the Lie algebra with
  non-trivial brackets $[v_1,v_2]=v_2+v_3$, $[v_1,v_3]=v_3$. Then
  $(S\frak{g})^{\frak{g}}_{\si}=k[v_3]$. On the other hand the kernel of the
weight of $v_3$ is the abelian Lie algebra $kv_2+kv_3$ whose semi-invariants
are $k[v_2,v_3]$. So this is different.  It turns out that in this case
we have to take $\frak{g}'=kv_1+kv_2+kv_3$ with non-trivial brackets
$[v_1,v_2]=v_3$.
\end{example}
The proof of Proposition \ref{ref-3.1-10} will be given after some preparation.
\begin{proposition} 
\label{ref-3.3-11} Assume that $\frak{h}$ is an ideal in $\frak{g}$ of
  codimension one.  Then one of the inclusions
$
(S\frak{h})^{\frak{g}}_{\si}\subset (S\frak{h})^{\frak{h}}_{\si}$ or 
$(S\frak{h})^{\frak{g}}_{\si}\subset (S\frak{g})^{\frak{g}}_{\si}$
is an equality. 
\end{proposition}
This result is perhaps better appreciated in the following equivalent formulation.
\begin{corollary} Assume that $\frak{h}$ is an ideal in $\frak{g}$ of
  codimension one.  Then either $(S\frak{h})^{\frak{h}}_{\si}\subset
  (S\frak{g})^{\frak{g}}_{\si}$ or  $(S\frak{g})^{\frak{g}}_{\si}\subset
  (S\frak{h})^{\frak{h}}_{\si}$
\end{corollary}
Note that \emph{both} inclusions may be equalities. This happens already for the two dimensional
non-abelian  Lie algebra.
\begin{proof}[Proof of Proposition \ref{ref-3.3-11}]
  Assume that the statement is false.  Write $\frak{g}$ as a
  semi-direct product $\frak{h}{+}kp$.  Let $f$ be a semi-invariant
  with weight $\chi$ in $(S\frak{h})^{\frak{h}}_{\si}\setminus
  (S\frak{h})^{\frak{g}}_{\si}$ which is irreducible in $S\frak{h}$
  and let $g=\sum_{i=0}^n a_i p^i$ be a semi-invariant with weight
  $\psi$ in $S\frak{g}$ such that $a_i\in S\frak{h}$, $n>0$ and
  $a_n\neq 0$.

  Since $f$ is not a semi-invariant for $\frak{g}$ it is not a
  semi-invariant for $p$.  From the fact that $g$ is a semi-invariant
  for $p$ we easily deduce that the $a_i$ are semi-invariants for $p$.
  In particular the non-zero $a_i$ cannot be divisible by $f$ (since a
  factor of a semi-invariant for $p$ is a semi-invariant for $p$).

We will now obtain a contradiction by computing the Poisson bracket
$\{f,g\}$.
Since $g$ is a semi-invariant in $S\frak{g}$ with weight $\psi$ we have
\[
\{f,g\}=\partial_\psi(f)\sum_i a_i p^i
\]
Since $f$ is a
semi-invariant in $S\frak{h}$ with weight $\chi$ we have
\[
\{f,g\}=-\sum_i \partial_\chi(a_i) fp^i+\sum_i a_i i p^{i-1} \{f,p\}
\]
Assume first $\partial_\psi(f)\neq 0$. Using the fact that
$\partial_\psi(f)$ has lower degree than $f$ and hence is not
divisible by $f$ we conclude that $a_n$ is divisible by $f$. Since
$a_n\neq 0$ this is a contradiction.

Assume now $\partial_\psi(f)=0$. In that case we obtain from the fact
that $f$ does not divide $\{f,p\}$ (since $f$ is not a semi-invariant
for $p$) that $f$ divides $a_i$ for $i>0$.  This is again a
contradiction.
\end{proof}
In the next two propositions we give some conditions under which the
ring of semi-invariants for a representation does not change under
passage to an ideal of the Lie algebra.
\begin{proposition} \label{ref-3.5-12} Let $V$ be a finite dimensional
representation of $\frak{g}$. Assume that $\frak{h}$ 
is an ideal in $\frak{g}$ such
that $r_{\frak{h}}(V) =r_{\frak{g}}(V)$. Then
$(SV)^{\frak{h}}_{\si}=(SV)^{\frak{g}}_{\si}$.
\end{proposition}
\begin{proof} 
  Assume $r_{\frak{h}}(V) =r_{\frak{g}}(V)$ and $(SV)^{\frak{h}}_{\si}\neq
  (SV)^{\frak{g}}_{\si}$.  Let $f\in
  (SV)^{\frak{h}}_{\si}\setminus (SV)^{\frak{g}}_{\si}$
  be a semi-invariant with weight $\chi$ which is irreducible in
  $SV$. Let $(h_i)_i$ be a basis
  of $\frak{h}$ and let $p\in \frak{g}-\frak{h}$ be such that 
$f$ is not a semi-invariant for $p$. 
Then by elementary
  linear algebra applied to the structure matrices of $V$ with respect
  to $\frak{g}$ and $\frak{h}$ there exist $a$, $b_i\in SV$
  with $a\neq 0$ such that $\delta=a\otimes p+\sum_i b_i\otimes h_i\in
  SV\otimes \frak{g}$ has the property that $\rho(\delta)$ acts
  trivially on $V$. Hence $\delta\in \ker \rho$.

  We claim we may choose $\delta$ in
  such a way that $a$ is not divisible by $f$.  Assume on the contrary
  that $a=f^na'$, $n>0$ such that $f$ does not divide $a'\in SV$.

 Since $\ker\rho$ is an ideal we
  have that $[1\otimes p,\delta]\in \ker\rho$ and
\[
[1\otimes p,\delta]=p(a)\otimes p+\sum_i b_i' \otimes h_i
\]
for suitable $b_i'\in SV$.  Then $p(a)=nf^{n-1}p(f)a'+f^np(a')$.
Since $p(f)$ is not divisible by $f$ (as it is not a semi-invariant for $p$) we see that the highest power of
$f$ which divides $p(a)$ is $f^{n-1}$ (and $p(a)\neq 0$). Replacing
$\delta$ by $[1\otimes p,\delta]$ and repeating this procedure we eventually
arrive at a $\delta$ such that $a$ is no longer divisible by $f$.

Applying this new $\delta$ to
$f$ we get
\[
0=ap(f)+\sum_i b_ih_i(f)=ap(f)+\sum b_i\chi(h_i)f
\]
Since neither $a$ nor $p(f)$ is divisible by $f$ we have obtained
a contradiction.
\end{proof}
\begin{proposition} \label{ref-3.6-13} Let $V$ be a finite dimensional representation of
  $\frak{g}$ and assume that $\frak{h}$ is an ideal in $\frak{g}$ of
  codimension one such that $\frak{g}=\frak{h}+ks$ with $s$ acting
  semi-simply on both $V$ and $\frak{h}$. Then
  $(SV)_{\si}^{\frak{h}}=(SV)_{\si}^{\frak{g}}$
\end{proposition}
\begin{proof} Put $S=SV$. We decompose $\frak{h}$ and $S$ according to
 the $s$-weights (i.e.\ as $s$-eigenspaces): $\frak{h}=\oplus_{\mu\in k}
  \frak{h}_\mu$, $S=\oplus_{\lambda\in k} S_\lambda$. For $f\in S$
let $\Supp f$ be the set of $\lambda$ such that $f_\lambda\neq 0$ in the decomposition
$f=\sum_{\lambda\in k}f_\lambda$ with $f_\lambda\in S_\lambda$. 

Let $f\in S$ be a semi-invariant for $\frak{h}$ with weight $\chi$. Write
$f=\sum_{\lambda\in k}f_\lambda$ $f_\lambda\in S_\lambda$. We claim
that the $f_\lambda$ are semi-invariants for $\frak{h}$, which implies
that they are in fact semi-invariants for $\frak{g}=\frak{h}{+}ks$. 
Hence $(SV)_{\si}^{\frak{h}}\subset (SV)_{\si}^{\frak{g}}$. Since
the other inclusion is obvious we are done. 

To prove the claim first assume $\chi=0$. Thus $f\in S^{\frak{h}}$.  Pick $h\in
  \frak{h}_\mu$. Then $0=h(f)=\sum_{\lambda\in k}h(f_\lambda)$ with
  $h(f_\lambda)\in S_{\mu+\lambda}$. Hence $h(f_\lambda)=0$ and thus $f_\lambda\in
S^{\frak{h}}$. So this case is OK. 

Now assume $\chi\neq 0$.  We first assert that $\chi(\frak{h}_\mu)=0$
for $\mu\neq 0$. To see this assume there exist $h\in \frak{h}_\mu$,
$\mu\neq 0$ such that $\chi(h)\neq 0$. From the equation
$h(f)=\chi(h)f$ we deduce $\Supp f=\Supp (h(f)) \subset \mu+\Supp
f$ which is impossible if $\mu\neq 0$.  So our assertion is correct.

As in the case $\chi=0$ we now deduce for $h\in \frak{h}_\mu$ that
$h(f_\lambda)=0$ if $\mu\neq 0$ and $h(f_\lambda)=\chi(h) f_\lambda$, if $\mu=0$.
So the $f_\lambda$'s are semi-invariants in this case also. 
\end{proof}
\begin{lemma} \label{ref-3.7-14} Assume that $f\in S\frak{g}$ is a
  semi-invariant with weight $\chi$ and $\frak{h}=\ker \chi$. Then
$c(\frak{g})=c(\frak{h})$.
\end{lemma}
\begin{proof} 
We may assume $\frak{g}\neq \frak{h}$, i.e. $\chi$ is non-trivial. 
Assume $c(\frak{g})\neq c(\frak{h})$.
Choose $p\in \frak{g}$ such that $\chi(p)=1$. Comparing
\begin{align*}
c(\frak{g})&=\dim \frak{g}-\frac{1}{2}r(\frak{g})\\
c(\frak{h})&=\dim \frak{g}-\frac{1}{2}r(\frak{h})-1
\end{align*}
we see that
\[
r(\frak{g})\neq r(\frak{h})+2
\]
Since the structure matrix of $\frak{g}$ is obtained from that of $\frak{h}$ by
adding a row and a column we have $r(\frak{g})-
r(\frak{h})\in \{0,2\}$. We obtain
\[
r(\frak{g})=r(\frak{h})
\]
The proof now parallels that of Lemma \ref{ref-3.5-12}. Let $(h_i)_i$ be a basis of $\frak{h}$ and select $a$, $b_i\in S\frak{g}$ with $a\neq 0$
such that $\delta=a\otimes p+\sum_i b_i\otimes h_i$ acts trivially on
$p$ and $(h_j)_j$. Thus $\delta\in \ker \rho$. 
In other words $\rho(\delta)$ acts trivially on $S\frak{g}$. 
But we also  find $\rho(\delta)(f)=af\neq 0$. This is a contradiction. 
\end{proof}
\begin{proof}[Proof of Proposition \ref{ref-3.1-10}]
  We will construct $\frak{g}'$ one step at a time. Assume that
  $\frak{g}$ has a non-trivial weight $\chi$ on $S\frak{g}$. Then we
  will construct a Lie algebra $\tilde{\frak{g}}$ such that $c({
    \tilde{\frak{g}} })=c(\frak{g})$,
  $(S\frak{g})_{\si}^{\frak{g}}=(S{ \tilde{\frak{g}} })_{\si}^{
    \tilde{\frak{g}} }$ and such that either $\dim {\tilde{\frak{g}}
  }<\dim \frak{g}$ or $\dim N( \tilde{\frak{g}} )>\dim N(\frak{g})$
  where $N(-)$ denotes the nil-radical. It is clear that by repeating
  this procedure we eventually end up with a Lie algebra which has the
  requested properties.

Let $f\in S\frak{g}$ be an non-zero eigenfunction for $\chi$. 
Put $\frak{h}=\ker \chi$. Since semi-invariants Poisson commute and since
the Poisson centralizer of $f$ is equal to $S\frak{h}$ we find (see \cite[Cor.\ 1.15]{Ooms3}).
\begin{equation}
\label{ref-3.1-15}
 (S\frak{g})_{\si}^{\frak{g}} \subset (S\frak{h})_{\si}^{\frak{h}}
\end{equation}
Choose $c\in \frak{g}$ such that $\chi(c)=1$. Then 
$
\ad(c)=D_s+D_p
$
where $D_s$, $D_p$ are two commuting derivations of $\frak{h}$ with $D_s$
being semi-simple and $D_p$ nilpotent. Let $\frak{j}=\frak{h}{+}ks{+}kp$ be
the semi-direct product of $\frak{h}$ with an abelian Lie algebra $ks{+}kp$
such that $\ad(s)$ acts by $D_s$ and $\ad(p)$ acts by $D_p$. Sending
$c$ to $s+p$ yields an embedding $\frak{g}\subset \frak{j}$. Put 
$\frak{k}=\frak{h}+kp$. Then we have $\frak{j}=\frak{g}+ks=\frak{k}+ks$.

Since $\ad(s)$ acts semi-simply on everything we have by Proposition \ref{ref-3.6-13}
\[
(S\frak{g})^{\frak{g}}_{\si}=(S\frak{g})^{\frak{j}}_{\si}=(S\frak{g})^{\frak{k}}_{\si}
\]
and thus by \eqref{ref-3.1-15}
\[
(S\frak{g})^{\frak{g}}_{\si}
=S\frak{h}\cap (S\frak{g})^{\frak{k}}_{\si}=
(S\frak{h})^{\frak{k}}_{\si}
\]
By Proposition \ref{ref-3.3-11} we have either $(S\frak{h})^{\frak{k}}_{\si}= (S\frak{h})^{\frak{h}}_{\si}$ or $(S\frak{h})^{\frak{k}}_{\si}=(S\frak{k})^{\frak{k}}_{\si}$.

Assume first $(S\frak{h})^{\frak{k}}_{\si}=
(S\frak{h})^{\frak{h}}_{\si}$. Then we put $\tilde{\frak{g}}=\frak{h}$. By
Lemma \ref{ref-3.7-14} we have $c(\tilde{\frak{g}})=c(\frak{g})$. Since
$\dim \frak{\tilde{g}}<\dim\frak{g}$ this case is done. 

Now assume $(S\frak{h})^{\frak{k}}_{\si}\neq
(S\frak{h})^{\frak{h}}_{\si}$ and thus $(S\frak{h})^{\frak{k}}_{\si}=(S\frak{k})^{\frak{k}}_{\si}$.  In this case we put
$\tilde{\frak{g}}=\frak{k}$ and hence we have $\dim \tilde{\frak{g}}= \dim
\frak{g}$. By Proposition \ref{ref-3.5-12} we have $r_{\frak{k}}
(\frak{h})>r_{\frak{h}}(\frak{h})$ and hence $r(\frak{k})=r_{\frak{k}}
(\frak{k})>r_{\frak{h}}(\frak{h})=r(\frak{h})$.

If $r(\frak{g})=r(\frak{h})$ then
$r_{\frak{h}}(\frak{g})=r_{\frak{h}}(\frak{h})$ and hence by
Proposition \ref{ref-3.5-12}
$(S\frak{h})^{\frak{g}}_{\si}=(S\frak{h})^{\frak{h}}_{\si}$. Since
$(S\frak{h})^{\frak{g}}_{\si}=
(S\frak{g})^{\frak{g}}_{\si}=
(S\frak{h})^{\frak{k}}_{\si}
$ 
this is a contradiction. Thus
$r(\frak{g})> r(\frak{h})$. Since the ranks involved jump at most by 2 we
deduce $r(\frak{g})=r(\frak{k})$ and hence $c(\tilde{\frak{g}})=c(\frak{g})$. 

It remains to show that in this case we have $\dim N(\frak{g})<\dim
N(\tilde{\frak{g}})$.  Since $p$ acts nilpotently we have
$N(\tilde{\frak{g}})=kp{+}N(\frak{h})$. We claim that
$N(\frak{g})\subset N(\frak{h})$ which is sufficient. To prove this
claim we need to show that no element of the form $c+n$ for $n\in \frak{h}$
acts nilpotently on $\frak{h}$.  Assume such $n$ exists. Then
$0=\chi(c+n)=\chi(c)$ and hence $c\in \frak{h}$ which is a contradiction.
\end{proof}
\section{A formula for the transcendence degree of invariants}
If $S$ is a commutative domain then we denote its field of
fractions by $Q(S)$.  In this section we prove the following result.
\begin{proposition}
\label{ref-4.1-16}
Let $V$ be a finite dimensional representation of $\frak{g}$ and assume
that $SV$ contains no proper semi-invariants. Then 
\begin{equation}
\label{ref-4.1-17}
\trdeg Q(SV)^{\frak{g}}=\dim V-r(V)
\end{equation}
\end{proposition}
In the case that $\frak{g}$ acts algebraically on $V$ the formula
\eqref{ref-4.1-17} was proved by Dixmier \cite[Lemme 7]{Dixmier} (it is a
more or less direct consequence of Rosenlicht's theorem
\cite{Rosenlicht}). Here we have traded algebraicity for the absence
of proper semi-invariants.  Both conditions are independent as Example
\ref{ref-4.8-25} below shows.

\medskip

Let $L/k$ be a finitely generated field extension on which $\frak{g}$ acts
by derivations. 
 Put $K=L^{\frak{g}}$. It is clear
that $K$ is algebraically closed in $L$.  We will
say that the action is \emph{geometric} if the induced map
$\rho:L\otimes_K \frak{g}\r \Der_K(L)$ is surjective.  If
$x_1,\ldots,x_n$ is a transcendence basis of $L/K$ and
$\partial_i=\partial/\partial x_i:L\r L$ are the corresponding
derivations then $\Der_K(L)=\sum_i L\partial_i$. From this we deduce
\begin{lemma} \label{ref-4.2-18} Let $L$ be as above.  Then the action is 
geometric if one has
\begin{equation}
\label{ref-4.2-19}
\trdeg_{L^{\frak{g}}} L=\dim\frak{g}-\dim_L \ker\rho
\end{equation}
\end{lemma}
Note that $\ker \rho$ can also be computed as
$\ker(L\otimes_k\frak{g}\r \Der_k(L))$.  So the number on the right hand
side of \eqref{ref-4.2-19} can be computed without knowing
$L^{\frak{g}}$.
\begin{lemma}
\label{ref-4.3-20}
Let $M\supset L\supset K$ be finitely generated field extensions of $k$
and let $\frak{g}$ be a finite dimensional Lie algebra acting on $M$
such that $K=M^{\frak{g}}$. Let $\frak{h}$ be an ideal in $\frak{g}$
and put $L=M^{\frak{h}}$.  Assume that the $\frak{h}$-action on $M$ is
geometric and likewise for the action $\frak{g}/\frak{h}$ on $L$. Then
the action of $\frak{g}$ on $M$ is geometric.
\end{lemma}
\begin{proof}
This follows from the following commutative diagram 
\[
\begin{CD}
0@>>> M\otimes_k \frak{h} @>>> M\otimes_k \frak{g}@>>> 
M\otimes_k \frak{g}/\frak{h} @>>> 0\\
@. @VVV @VVV @VVV @.
\\
0@>>>\Der_L(M) @>>> \Der_K(M) @>>> M\otimes_L \Der_K(L) @>>> 0\\
@. @VVV @. @VVV @.\\
@. 0 @. @. 0 
\end{CD}
\]
\end{proof}
Now let $S$ be a $k$-algebra which is an integral domain
with finitely generated fraction field. 
Assume $\frak{g}$ acts on $S$. We
say that $\frak{g}$ acts \emph{generically geometrically}  if the induced
action on the fraction field is geometric. 

If $Y$ is a (possibly singular) variety then we denote the tangent
space in a point $y\in Y$ by $T_{Y,y}$. If $Y$ is smooth then $T_{Y,y}$ is
the fiber of the tangent bundle $T_Y$ of $Y$ at $y$.
\begin{proposition} 
  \label{ref-4.4-21} Assume that $S/k$ is a finitely generated
  domain and that $\frak{g}$ acts on $S$. Put $Y=\Spec S$. Let $L$ be
  the fraction field of $S$. Then the action is generically geometric
  if and only if
\[
\trdeg_{L^{\frak{g}}} L=\dim\frak{g}-\min_{y\in Y}\dim \frak{g}_y
\]
where the minimum is taken over the closed points in $Y$ and $\frak{g}_y$
denotes the stabilizer of $y\in Y$, i.e.\ the kernel of $\frak{g}\r T_{Y,y}$.
\end{proposition}
\begin{proof} For technical reasons it is more convenient to work
  with differentials instead of with vector fields as the sheaf of
  differentials is always compatible with taking fibers.  

There is a
  canonical pairing $\frak{g}\otimes_k\Omega_Y\r \Oscr_Y: v\otimes
  f\,dg\mapsto fv(g)$ which yields a map of coherent $\Oscr_Y$ modules
  $\rho^\ast:\Omega_Y\r \Oscr_Y\otimes \frak{g}^\ast$. Taking the
  fiber in a point $y\in Y$ one checks $\frak{g}^\ast_y=\coker
  \rho^\ast_y$. By semi-continuity the dimension of the cokernel of
  the generic fiber of $\rho^\ast$ is equal to the minimum of the
  dimensions of the cokernels of the special fibers.

Thus we find
\begin{align*}
\dim_L \ker(L\otimes_k\frak{g}\r \Der_k(L))&=\dim_L \coker(\Omega_{L,k}
\r L\otimes_k\frak{g}^\ast)\\&=
\min_{y\in Y}\dim \frak{g}_y^\ast\\ 
&=
\min_{y\in Y}\dim \frak{g}_y
\end{align*}
It now suffices to apply Lemma \ref{ref-4.2-18}. 
\end{proof}
\begin{lemma} \label{ref-4.5-22} Assume that $\frak{g}=\Lie(G)$ where
  $G$ is a connected algebraic group acting rationally on a domain $S$
with finitely generated fraction field. Then
  the action is generically geometric.
\end{lemma}
\begin{proof} Let $L$ be the fraction field of
  $S$. Since $L$ is finitely generated and since $G$ acts rationally
  on $S$ we may select a finitely generated $G$-invariant subring
  $S_0\subset S$ such that $L$ is the fraction field of $S_0$.  From here on
the proof proceeds as in \cite[Lemme 7]{Dixmier1}. For the benefit of
the reader let us repeat the argument. By
  Rosenlicht's theorem there exists an $G$-invariant open $U\subset
  Y=\Spec S_0$ such that a geometric quotient $U/G$ exists. Shrinking
  $U$ we may assume that $U$ is smooth. By the properties of a geometric quotient the field of rational
  functions on $U$ is $L$, the field of rational functions on $U/G$ is
  $L^G$ and the fibers of $U\r U/G$ are the $G$-orbits. Hence $\trdeg
  L/L^G=\dim U-\dim U/G$ is equal to the dimension of the generic $G$
  orbit on $U$ or equivalently on $Y$. This is $\dim
  G-\min_{y\in Y} \dim G_y =\dim \frak{g}-\min_{y\in Y}\dim\frak{g}_y$. We may
now apply Proposition \ref{ref-4.4-21} with $S=S_0$.
\end{proof}
We would like to have a
 transitivity result as in Lemma \ref{ref-4.3-20} but one does not
always have $Q(S)^{\frak{g}}=Q(S^{\frak{g}})$. So we introduce a special
situation in which this identity holds.

Let us say that a graded ring $S$ is \emph{connected} if it is of the
form $k+S_1+S_2+\cdots$ with $\dim S_i<\infty$. For technical reasons
we do not assume that $S$ is finitely generated. 
Recall the following.
\begin{lemma} \label{ref-4.6-23} Let $S$ be a connected graded factorial
  domain, and assume that $\frak{g}$ acts in a graded way on $S$
  without proper semi-invariants. Then $S^{\frak{g}}$ is factorial,
  and furthermore $Q(S)^{\frak{g}}=Q(S^{\frak{g}})$.
\end{lemma}
The following is a weak version of Rosenlicht's theorem which is
also valid for non-algebraic Lie algebras.
\begin{proposition} 
\label{ref-4.7-24} Let $S$ be a connected factorial graded domain
  with finitely generated quotient field and assume that $\frak{g}$
  acts in a graded way on $S$ without proper
  semi-invariants. Then the action is generically geometric.
\end{proposition}
\begin{proof} We have a filtration by ideals
\[
0=\frak{g}_0\subset \cdots \subset \frak{g}_{n-1}\subset \frak{g}_n=\frak{g}
\]
where $\frak{g}_n/\frak{g}_{n-1}$ is semi-simple, and the other quotients
are abelian. We have a corresponding filtration
\[
S=S^0\supset \cdots \supset S^{n-1}\supset S^n=S^{\frak{g}}
\]
with $S^i=S^{\frak{g}_i}$ and hence $S^{i+1}=(S^{i})^{\frak{g}_{i+1}/\frak{g}_i}$. 

By Lemmas \ref{ref-4.6-23} and \ref{ref-4.3-20} we may assume that
$\frak{g}$ is either semi-simple or abelian. If $\frak{g}$ is
semi-simple then it acts algebraically (since $\dim S_i<\infty$) and
hence we may invoke Lemma \ref{ref-4.5-22}. If $\frak{g}$ is abelian
then its generalized weights must be zero (for otherwise we could construct
a proper semi-invariant). Hence $\frak{g}$ acts locally nilpotently
and hence algebraically. We may again invoke Lemma \ref{ref-4.5-22}.
\end{proof} 
\begin{proof}[Proof of Proposition \ref{ref-4.1-16}]
  Put $L=Q(SV)$.  By \eqref{ref-4.2-19} $\trdeg_{L^{\frak{g}}} L=\dim
  \frak{g} -n$ where $n$ is the dimension of the null space of the
  structure matrix. The structure matrix has size $\dim
  V\times \dim \frak{g}$. Thus $\trdeg_{L^{\frak{g}}} L=r(V)$.  Hence
  $\trdeg L^{\frak{g}}=\dim V-r(V)$.
\end{proof}
\begin{example} 
\label{ref-4.8-25} Let $W$ be a finite dimensional vector space.  Then
  the \emph{Heisenberg Lie algebra} $\frak{h}$ on $W$ is the vector
  space $W\oplus W^\ast\oplus kc$ with for all $w,w'\in W$, $\phi,\phi'\in
  W^\ast$: $[w,w']=0$, $[\phi,\phi']=0$, $[c,w]=0$, $[c,\phi]=0$,
  $[\phi,w]=\phi(w)c$.

  For $p\in \End(W)$ let $D_p$ be the derivation $(p,-p^\ast,0)$ of
  $\frak{h}$ and let $\frak{g}=\frak{h}{+}kt$ be the corresponding
  semi-direct product. Put $S=S\frak{g}$.  Assume that $p$ is
  invertible. Then $\frak{g}$ has a non-degenerate invariant symmetric
  bilinear form $(-,-)$ whose non-trivial values are given by
  $(t,c)=1$, $(\phi,w)=-\phi(p^{-1}(w))$. Hence $\frak{g}$ is
  quadratic. It then follows from \cite[Cor.\ 2.3, Prop.\ 3.2]{Ooms4}
  that $S\frak{g}$ contains no proper semi-invariants.  However if we
  choose $p$ to be non-diagonalizable then $\frak{g}$ does not act
  algebraically as $t$ does not have a Jordan decomposition.  So the
  hypotheses of Proposition \ref{ref-4.1-16} do not imply that
  $\frak{g}$ acts algebraically.

  The formula \eqref{ref-4.1-17} asserts that the transcendence degree
  of $Q(S)^{\frak{g}}$ should be $2$, as an
  easy computation shows $r(\frak{g})=\dim \frak{g}-2$.

  This example is simple enough to verify directly. As an aside we find
  that the hypothesis that $p$ is invertible is in fact superfluous.
  The absence of proper semi-invariants always holds. 

 Choose a basis
  $(w_i)_i$ for $W$ and assume that $p(w_i)=\sum_{j}p_{ij}w_j$. Then
  one has $p^\ast(w_j^\ast)=\sum_i p_{ij} w_i^\ast$. Put
  $z=tc-\sum_{ij} p_{ij}w_i w_j^\ast$ (this is the Casimir element for the pairing $(-,-)$ in the 
case that $p$ is invertible).  Then $c,z\in S^{\frak{g}}$.
  We write $S_c=k[c^{\pm 1},z,(w_ic^{-1})_i,w_j^{\ast}]$. With respect
  to these new generators the only non-trivial Poisson brackets are
  $\{w_j^\ast,w_ic^{-1}\}=\delta_{ji}$. A semi-invariant in $S$
  generates a principal Poisson ideal in $S$ and hence in $S_c$.
  Computing with the new generators we see that the only principal
  Poisson ideals in $S_c$ are those generated by elements in $k[c^{\pm
    1},z]$. Thus $S^{\frak{g}}_{\si}=S\cap k[c^{\pm 1},z]$ which is
  equal to $k[c,z]$ if $p\neq 0$ and equal to $k[c,t]$ otherwise.
  Thus we see that $S\frak{g}$ contains no proper semi-invariants.
In both cases we find  $\trdeg Q(S)^{\frak{g}}=2$ as predicted by \eqref{ref-4.1-17}.
 \end{example}
\section{Proofs in the absence of proper semi-invariants}
Throughout $V$ is a finite dimensional representation of $\frak{g}$.
For convenience we write $S=SV$, $R=(SV)^{\frak{g}}$ and we let $L$ be
the field of fractions of $S$.  As we will mostly use geometrical
language we also put $Y=\Spec S=V^\ast$, $X=\Spec R$ and we let
$\pi:Y\r X$ be dual to the inclusion $R\r S$. If $R$ is finitely generated then the regular locus of
$X$ is denoted\footnote{Unfortunately the subscript ``$\text{reg}$'' is already taken
by $\frak{g}^\ast_{\text{reg}}$.} by $X_{\sm}$. 

The following result is an adaptation of \cite{Knop3} to the case of 
non-semisimple Lie algebras. 
\begin{proposition} 
\label{ref-5.1-26} 
Assume that $(SV)_{\si}^{\frak{g}}=(SV)^{\frak{g}}$ and that $(SV)^{\frak{g}}$ is
finitely generated.
Let 
\[
U=\{y\in Y\mid \text{$\pi(y)\in X_{\sm}$ and $\pi$
    is smooth in $y$}\}
\]
Then $\codim_Y(Y-U)\ge 2$.
\end{proposition}
\begin{proof}
  The proof is that of \cite{Knop3} with minor adaptations.  Without
  loss of generality we may replace $\frak{g}$ by the algebraic hull
  of the image of $\frak{g}$ in $\End_k(V)$.  Let $G\subset \GL(V)$ be
  the affine connected algebraic group such that $\Lie G=\frak{g}$.
  Then $G$ acts rationally on $S$ and $R=S^G$.

Let $E$ be the union of the irreducible divisors in $Y-U$. Then $E$ is
$G$-invariant. Since $G$ is connected it follows that $E$ is
irreducible. Since $S$ is factorial it follows that $E=V(f)$ for some
irreducible $f\in S$.

For $\sigma \in G$ we
have that $c_\sigma=\sigma(f)f^{-1}$ is a unit in $S$ and hence 
$c_\sigma\in k^\ast$.
Thus $f$ is a semi-invariant and hence $f\in R$.

We claim that the map
$
R/fR\r S/f S
$
is injective. Assume there is some element $\bar{c}$ in the kernel. Then
$c=fd$ with $d\in S$. But then $d\in S^{\frak{g}}=R$. Hence $\bar{c}=0$.

Let $D$ be the divisor in $X$ of $f$. Then $E=\pi^{-1}(D)$ and the map
$E\r D$ is a dominant map between irreducible algebraic varieties. Since
$X$ is normal \hbox{$D\cap X_{\sm}\neq \emptyset$}. 

 By
generic smoothness there exist dense open $E'\subset E$, $D'\subset D\cap X_{\sm}$
such that $E',D'$ are regular and $\pi$ restricts to a smooth map
$E'\r D'$.

Let $y\in E'$. We will show that $\pi$ is smooth at $y$, contradicting the 
fact that $E'$ is contained in the non smooth locus of $\pi$. We consider
the following commutative diagram of tangent spaces with $x=\pi(y)$
\[
\begin{CD}
0 @>>> T_{E',y} @>>> T_{Y,y} @>df_y>> k @>>> 0\\
@. @V d\pi\mid_{E'} VV @V d\pi VV @| \\
0 @>>> T_{D',x} @>>> T_{X,x} @>>df_x> k @>>> 0
\end{CD}
\]
Since $x$ is regular in $D$, $X$ and $y$ is regular in $E$, $Y$  the rows are exact.  The left
most map is surjective since $E'\r D'$ is smooth. This implies that
the middle map is surjective.
\end{proof}
We keep the notations as in the statement of Proposition \ref{ref-5.1-26}
and we assume throughout that $(SV)^{\frak{g}}_{\si}=(SV)^{\frak{g}}$ and that 
$(SV)^{\frak{g}}$ is finitely generated. 
This implies in particular that $(SV)^{\frak{g}}$ is factorial and
$L^{\frak{g}}=Q((SV)^{\frak{g}})$. Furthermore by Propositions
\ref{ref-4.7-24} and \ref{ref-4.4-21} we have
\begin{equation}
\label{ref-5.1-27}
\dim Y-\dim X=\trdeg_{L^{\frak{g}}} L=\dim\frak{g}-\min_{y\in Y}\dim \frak{g}_y
\end{equation}
This leads to the definition
\[
Y'=\{y\in Y\mid \dim \frak{g}-\dim \frak{g}_y<\dim Y-\dim X\}\subsetneq Y
\]
and we will also put
\begin{equation}
\label{ref-5.2-28}
W=U\cap (Y-Y')
\end{equation}
Since the elements of $\frak{g}$ define vector fields on $Y$ which
annihilate invariant functions, we have the usual  map of vector bundles on $Y$
\[
\rho:\Oscr_Y \otimes_k \frak{g}\r T_Y
\] 
which extends to a complex of vector bundles on~$\pi^{-1}X_{\sm}$
\begin{equation}
\label{ref-5.3-29}
 \Oscr_{\pi^{-1}X_{\sm}} \otimes_k \frak{g}
\xrightarrow{\rho_{\pi^{-1}X_{\sm}}} T_{\pi^{-1}X_{\sm}}
\xrightarrow{d\pi_{\pi^{-1}X_{\sm}}} \pi^\ast_{\pi^{-1}X_{\sm}} T_{X_{\sm}}
\r 0
\end{equation}
We see that $U$  is the locus in $\pi^{-1}X_{\sm}$ where \eqref{ref-5.3-29} 
is exact at $\pi^\ast_{\pi^{-1}X_{\sm}} T_{X_{\sm}}$ and $W$ is the locus where the entire
complex is exact. Thus we have an exact sequence
\begin{equation}
\label{ref-5.4-30}
\Oscr_W \otimes_k  \frak{g}  \xrightarrow{\rho_W} T_W\xrightarrow{d\pi_W} \pi^\ast_W T_{X_{\sm}}
\r 0
\end{equation}
The following is one of the main results of \cite{Knop3}. For completeness
we include the proof in our setting. 
\begin{proposition} \label{ref-5.2-31} One has
$W=(Y-Y')\cap \pi^{-1}X_{\sm}$. 
\end{proposition}
The statement of this Proposition means that if $y\in Y$ is such that
$\pi(y)=x$ is regular in $X$ and $\frak{g}_y$ has minimal dimension
then $\pi$ is smooth in $y$.
\begin{proof}[Proof of Proposition \ref{ref-5.2-31}]
  Put $\overline{W}=(Y-Y')\cap \pi^{-1}X_{\sm}$. Since $W=U\cap
  (Y-Y')$ the inclusion $W\subset \overline{W}$
  is obvious. To prove the opposite inclusion we look at the complex
\[
\Oscr_{\overline{W}} \otimes_k{\frak{g}} \xrightarrow{\rho_{\overline{W}}} T_{\overline{W}}  \xrightarrow{d\pi_{\overline{W}}} \pi^\ast_{\overline{W}} T_{X_{\sm}}
\r 0
\]
Since $\rho_{\overline{W}}$ has constant rank $\coker \rho_{\overline{W}}$ is a
vector bundle. Furthermore by the exactness of \eqref{ref-5.4-30} we deduce
that $\coker \rho_{\overline{W}}\r \pi^\ast_{\overline{W}} T_{X_{\sm}}$ is an isomorphism on $W$.

Since $W\subset \overline{W}\subset Y-Y'$ and $Y-Y'-W=(Y-Y')\cap (Y-U)$
has codimension $\ge 2$ in $Y-Y'$, the same holds for the
codimension of $\overline{W}-W$ in $\overline{W}$.  It follows that
$\coker \rho_{\overline{W}}\r \pi^\ast_{\overline{W}}
T_{X_{\sm}}$ is an isomorphism on the whole of $\overline{W}$. In particular the map
$T_{\overline{W}}\r \pi^\ast_{\overline{W}} T_{X_{\sm}}$ is
surjective.  This implies that $\overline{W}\subset U$. Hence
$\overline{W}=W$.
\end{proof}

\begin{lemma} (see also \cite[Lemma 4]{Knop1}) \label{ref-5.3-32}
  The dual of \eqref{ref-5.3-29} yields an exact sequence
\begin{equation}
\label{ref-5.5-33}
0\r \pi^\ast_{\pi^{-1}X_{\sm}}\Omega_{X_{\sm}}\xrightarrow{d\pi^\ast_{\pi^{-1}X_{\sm}}}
 \Omega_{\pi^{-1}X_{\sm}}\xrightarrow{\rho^\ast_{\pi^{-1}X_{\sm}}} \frak{g}^\ast\otimes \Oscr_{\pi^{-1}X_{\sm}}
\end{equation}
of vector bundles on $\pi^{-1}X_{\sm}$. 
\end{lemma}
\begin{proof}
 Let $C=\coker \rho_{\pi^{-1}X_{\sm}}$. Then we have commutative diagram
\[
\begin{CD}
\Oscr_{\pi^{-1}X_{\sm}} \otimes_k  \frak{g}
@>\rho_{\pi^{-1}X_{\sm}}>> T_{\pi^{-1}X_{\sm}} @>d\pi_{\pi^{-1}X_{\sm}}>> \pi^\ast_{\pi^{-1}X_{\sm}} T_{X_{\sm}}
\\
@| @| @AA\beta A\\
\Oscr_{\pi^{-1}X_{\sm}} \otimes_k \frak{g}
@>>\rho_{\pi^{-1}X_{\sm}}> T_{\pi^{-1}X_{\sm}} @>>> C @>>> 0\\
@. @. @AAA\\
@. @. K\\
@. @. @AAA\\
@. @. 0
\end{CD}
\]
Since $d\pi_{\pi^{-1}X_{\sm}}$ is surjective on $U$, the same holds for $\beta$. 
So the vertical sequence is exact on~$U$.
Furthermore since the upper sequence is exact on $W$ it follows that
$K$ is torsion.  Dualizing the vertical sequence we obtain $C^\ast|U\cong
(\pi^\ast_{\pi^{-1}X_{\sm}}T_{X_{\sm}})^\ast|U$. This extends to an isomorphism 
$C^\ast= (\pi^\ast_{\pi^{-1}X_{\sm}}  T_{X_{\sm}})^\ast$.

The lemma now follows by
dualizing the lower exact sequence.
\end{proof}
We can now prove a generalization of Proposition \ref{ref-1.6-6}(1)
(taking into account that if $V=\frak{g}$ we have $\rho^\ast=-\rho$, see
below). 
\begin{proposition} 
\label{ref-5.4-34}
  Assume that $(SV)_{\si}^{\frak{g}}=(SV)^{\frak{g}}$,
and that  $(SV)^{\frak{g}}$ is a (necessarily finitely generated) polynomial
  ring. Then $ \ker \rho^\ast $ is free.
\end{proposition}
\begin{proof} This follows immediately from Lemma \ref{ref-5.3-32}, taking
into account that $X=X_{\sm}$.
\end{proof}
Now we specialize to the case $V=\frak{g}$ but we still assume that
the conditions of Proposition \ref{ref-5.1-26} hold. Our assumption
that $(S\frak{g})_{\si}^{\frak{g}}=(S\frak{g})^{\frak{g}}$ implies in particular
that $\frak{g}$ is unimodular by \cite{DDV}. To lighten the
notations we silently fix an isomorphism $\wedge^{\dim
  \frak{g}}\frak{g}\cong k$.

In what follows
we need to keep track of the grading. Therefore we introduce an
additional 1-dimensional torus $\GG_m$ which acts with weight $n$ on the degree $n$-part of $S$.  Everything we do is equivariant with respect to this torus.
If  $L$ is the one dimensional representation of $\GG_m$
corresponding to the identity character
then we write $?(n)$ for $?\otimes L^{n}$. 

Taking into account $V=\frak{g}$ we obtain $T_Y= (\Oscr_Y \otimes
\frak{g}^\ast)(1)$, $\Omega_Y=(\Oscr_Y \otimes \frak{g})(-1)$. 
In particular $\rho$ may be viewed as a map:
\[
\rho: \Oscr_Y \otimes \frak{g} \r (\Oscr_Y \otimes  \frak{g}^\ast )(1)
\]
Since $\rho$ is represented by the structure matrix of
$\frak{g}$ which is anti-symmetric we obtain $\rho^\ast=-\rho(-1)$.
Concatenating \eqref{ref-5.3-29} with \eqref{ref-5.5-33} we obtain a complex
\begin{equation}
\label{ref-5.6-35}
0\r (\pi^\ast_{\pi^{-1}X_{\sm}}\Omega_{X_{\sm}})(1) \xrightarrow{d\pi^\ast}
\Oscr_{\pi^{-1}X_{\sm}} \otimes_k\frak{g}
\xrightarrow{\rho} 
( \Oscr_{\pi^{-1}X_{\sm}}\otimes \frak{g}^\ast
)(1)
\xrightarrow{d\pi} \pi^\ast_{\pi^{-1}X_{\sm}} T_{X_{\sm}}
\r 0
\end{equation}
 which is possibly non-exact at $(\Oscr_{\pi^{-1}X_{\sm}} \otimes \frak{g}^\ast )(1)$
and $\pi^{\ast}_{\pi^{-1}X_{\sm}}T_{X_{\sm}}$.  The locus where the right non-trivial map is surjective is $U$.
The locus
where it is exact is precisely~$W$. 

If $Z$ is a normal variety and $F$ is a coherent torsion
free $\Oscr_Z$-module then we put $\det F=(\wedge^{\rk F} F)^{\ast\ast}$. The
operation $F\mapsto \det F$ is multiplicative on complexes which are
exact in codimension $\ge 2$.  Furthermore the following is well
known.
\begin{lemma} 
\label{ref-5.5-36}
Let $\rho:F\r G$ be a map between vector bundles on $Z$ which is 
  generically of rank $r>0$. Let $M=\im \rho$ and let $N\subset G$ be
  the maximal coherent subsheaf of $G$ containing $M$ such that $N/M$ is
  torsion. Finally let $I(\rho)$ be the ideal in $\Oscr_Z$ locally
generated by the $r\times r$ minors in a matrix representation
of $\rho$.  Then $\det M=(I(\rho) \det N)^{\ast\ast}$ as submodules of
$\wedge^r G$. 
\end{lemma}
\begin{proof} We may reduce to the case that $Z$ is the spectrum of a
  discrete valuation ring $D$, $F=D^p$, $G=D^q$. In that case $\rho$
  may be diagonalized as
\[
\begin{pmatrix}
\pi^{a_1} &          0 & \cdots  &  0        &   0    &\cdots  & 0\\
0        &   \pi^{a_2} & \cdots  &  0        &   0    & \cdots & 0\\
\vdots   &    \vdots  &         &  \vdots   & \vdots &         & \vdots \\
0        &    0       & \cdots  &  \pi^{a_r} &   0     &\cdots   & 0\\
0        &    0       & \cdots  &  0        &   0     &\cdots   & 0\\
\vdots   & \vdots     &         & \vdots    & \vdots &          &\vdots\\
0        &    0       & \cdots  &  0        &   0     &\cdots   & 0
\end{pmatrix}
\]
where $\pi$ is a uniformizing element of $D$.  Thus
$I(\rho)=(\pi^{\sum_i a_i})$.  We find
$M=\pi^{a_1}D\oplus\cdots\oplus \pi^{a_r}D\oplus 0\oplus \cdots \oplus
0\subset D^q$ and $N=D^r\oplus 0\oplus\cdots \oplus 0\subset D^q$.
Let $(e_i)_i$ be the standard basis of $D^q$.
Then 
$\det M=\pi^{\sum_i a_i} D e_1\wedge\cdots \wedge e_r $, $\det N=De_1\wedge\cdots \wedge e_r $ and so we have indeed $\det M=I(\rho)\det N$ inside
$\wedge^r G$.  
\end{proof}
We put $ \omega_Z=\det \Omega_Z $.  This is the so-called
\emph{dualizing module} on $Z$. 

\medskip

Put $M=\im \rho$, $N=\ker d\pi$ in \eqref{ref-5.6-35}. We obtain two
exact sequences of torsion free $\Oscr_{\pi^{-1}X_{\sm}}$-modules. 
\[
0\r (\pi^\ast_{\pi^{-1}X_{\sm}}\Omega_{X_{\sm}})(1)\xrightarrow{d\pi^\ast}
 \Oscr_{\pi^{-1}X_{\sm}}\otimes_k\frak{g}
\r M\r 0
\]
\[
0\r N \r (\Oscr_{\pi^{-1}X_{\sm}}\otimes \frak{g}^\ast
)(1)
\xrightarrow{d\pi} \pi^\ast_{\pi^{-1}X_{\sm}} T_{X_{\sm}}
\]
where $\coker d\pi$ is supported on $\pi^{-1}X_{\sm}-U$ which has
codimension $\ge 2$ by Proposition \ref{ref-5.1-26}.  We have
$M\subset N$ and furthermore the support of $N/M$ is
$\pi^{-1}X_{\sm}-W$. Since $\pi^\ast_{\pi^{-1}X_{\sm}} T_{X_{\sm}}$ is
torsion free we find that $N$ is the maximal submodule of
$(\Oscr_{\pi^{-1}X_{\sm}}\otimes \frak{g}^\ast)(1)$ containing $M$
such that $N/M$ is torsion. 
Hence we are in the setting of Lemma
\ref{ref-5.5-36}, provided $\frak{g}$ is non-abelian (we want
$\rk \rho>0$), which we will
temporarily assume.  Let $f$ be the fundamental semi-invariant (cfr
\S\ref{ref-2-8}) in $S\frak{g}$. Then we have $I(\rho)\subset (f)$ and
$(f)/I(\rho)$ is supported in codimension $\ge 2$. Hence in the
application of Lemma \ref{ref-5.5-36} we may replace $I(\rho)$ by
$(f)=\Oscr_Y(-d(\frak{g}))$.

Taking into account $\dim \frak{g}=\dim Y$ we
conclude (using multiplicativity of $\det$)
\begin{align*}
\det M&=     (\pi^\ast_{\pi^{-1}X_{\sm}} \omega_{X_{\sm}})^\ast(-\dim X)\\
\det N&=(\pi^\ast_{\pi^{-1}X_{\sm}} \omega_{X_{\sm}})(\dim Y)
\end{align*}
Hence by Lemma \ref{ref-5.5-36} we obtain
\[
(\pi^\ast_{\pi^{-1}X_{\sm}} \omega_{X_{\sm}})^\ast(-\dim X)=
(\pi^\ast_{\pi^{-1}X_{\sm}} \omega_{X_{\sm}})(\dim Y-d(\frak{g}))
\]
or in other words
\begin{equation}
\label{ref-5.7-37}
\Oscr_{\pi^{-1}X_{\sm}}(-\dim Y-\dim X+d(\frak{g}))=\pi^\ast_{\pi^{-1}X_{\sm}}\omega_{X_{\sm}}^{\otimes 2}
\end{equation}
If $T=k+T_1+\cdots+$ is a finitely generated positively graded
normal commutative ring and if $\omega_T\cong T(-a)$ then we will call $a$ the
\emph{Gorenstein invariant} of $T$ and denote it by $a(T)$.  This is for
example always defined if $T$ is factorial.\footnote{This is a  slight
  abusing of existing terminology as normally the Gorenstein
  invariant is only defined for Gorenstein rings.}
\begin{example} \label{ref-5.6-38} \begin{enumerate}
\item If $T=k[f_1,\ldots,f_r]$ is a graded polynomial ring with homogeneous
generators $f_1,\ldots,f_r$ of strictly positive degree then 
$a(T)=\sum_i \deg f_i$.
\item Similarly if $T=k[f_1,\ldots,f_r]/(p_1,\ldots,p_s)$ is a
  homogeneous normal complete intersection then $a(T)=\sum_i \deg
  f_i-\sum_j\deg p_j$.
\end{enumerate}
\end{example}
We can now prove a more general version of Proposition \ref{ref-1.4-3} in the
absence of proper semi-invariants. 
\begin{proposition} \label{ref-5.7-39} Assume that 
  $(S\frak{g})^{\frak{g}}_{\si}= (S\frak{g})^{\frak{g}}$,
and   $(S\frak{g})^{\frak{g}}$ is finitely generated. Then
$a((S\frak{g})^{\frak{g}})$ is defined and is equal to
\begin{equation}
\label{ref-5.8-40}
a((S\frak{g})^{\frak{g}})=\frac{1}{2}(\dim \frak{g}+i(\frak{g})-d(\frak{g}))
\end{equation}
\end{proposition}
\begin{proof}
  We use the same notations as above. We assume that $\frak{g}$ is non-abelian
since otherwise the result is trivial. 
Since there are no proper
  semi-invariants, $(S\frak{g})^{\frak{g}}$ is factorial, and hence
  $a((S\frak{g})^{\frak{g}})$ is defined. Put $a=a((S\frak{g})^{\frak{g}})$.
Thus $\omega_X=\Oscr_X(-a)$. 

 Let $i:\pi^{-1}X_{\sm}\r
Y$ be the inclusion map.  Applying $i_\at$ to \eqref{ref-5.7-37} and using
the fact that by Proposition \ref{ref-5.1-26} $\codim_Y(Y-\pi^{-1}X_{\sm})\ge 2$ and that everything
is reflexive, we obtain an equality
\[
\Oscr_{Y}(-\dim Y-\dim X+d(\frak{g}))=\pi^\ast
\omega_{X}^{\otimes 2}=\Oscr_Y(-2a)
\]
and hence $2a=\dim X+\dim Y-d(\frak{g})=i(\frak{g})+\dim \frak{g}-d(\frak{g})$
which yields \eqref{ref-5.8-40}.
\end{proof}
This result can also be proved using the method exhibited in
\cite[Remark 1.6.3]{Panyushev} as Proposition \ref{ref-5.1-26} shows
that the set $S$ in loc.\ cit. (which is $U$ in our terminology) is
``big'' in the sense of \cite{Panyushev}.
\begin{example}
\label{ref-5.8-41}
  We apply Proposition \ref{ref-5.7-39} to a
  non-coregular example. Let $\frak{g}=L(6)$ (cfr Example \ref{ref-1.7-7}).
  Using \cite[\S89,\S93]{GY} or the library ``ainvar.lib'' from Singular
\cite{GPS05} we find $(S\frak{g})^{\frak{g}}=k[f_1,f_2,f_3,f_4,f_5,]/(p)$
where
\begin{align*}
f_1 &= v_5^2 - 2v_4v_6\\
f_2 &= v_5^3 - 3v_4v_5v_6 + 3v_3v_6^2\\
f_3 &= v_4^2 - 2v_3v_5 + 2v_2v_6\\
f_4 &= 2v_4^3 + 6v_2v_5^2  + 9v_3^2v_6  - 12v_2v_4v_6 - 6v_3v_4v_5\\
f_5 &= v_6
\end{align*}
and
\[
p=f_4f_5^3 -3f_1f_3f_5^2 +f_1^3 -f_2^2
\]
According to Example \ref{ref-1.7-7} we have $\dim \frak{g}=6$, 
$i(\frak{g})=4$, $d(\frak{g})=0$. The equality \eqref{ref-5.8-40} becomes
(taking into account Example \ref{ref-1.7-7}(2))
\[
2 + 3 + 2 + 3+1-6=5=(6+4-0)/2
\]
\end{example}
\begin{remark}
  Let $g^2$ be the fundamental semi-invariant in $S\frak{g}$ (it is a square
by Lemma \ref{ref-2.1-9}). If
  $(S\frak{g})_{\si}^{\frak{g}}$ is a polynomial algebra then it seems
  that in many cases the irreducible factors of $g$ form a subset of
  the generators of $(S\frak{g})_{\si}^{\frak{g}}$.  This is true for
  Frobenius Lie algebras \cite{Ooms3} and also for the examples
  covered by the methods in \cite{Ooms2}.

Assume that $(S\frak{g})^{\frak{g}}_{\si}=(S\frak{g})^{\frak{g}}$ and
$(S\frak{g})^{\frak{g}}=k[f_1,\ldots,f_r]$. If $g=\prod_i f_i^{\epsilon_i}$ then
the equality \eqref{ref-5.8-40} becomes 
\[
\sum_{i=1}^r (1+\epsilon_i) \deg f_i=c(\frak{g})
\]
This is similar to a phenomenon observed by Fauquant-Millet and Joseph 
in \cite{JFM}
that one can sometimes make the inequality  \eqref{ref-1.1-1} into an equality
by changing the degrees of the $f_i$ in some natural way. 
\end{remark}
We end this section by proving Proposition \ref{ref-1.6-6}(2)(3).
\begin{proposition}
\label{ref-5.10-42}
Assume that $(S\frak{g})^{\mathfrak{g}}_{\si}=(S\frak{g})^{\mathfrak{g}}$,
$\frak{g}$ is not abelian
  and $\frak{g}$ is coregular. Then $
\codim_{\frak{g}^\ast}(\frak{g}^\ast-\frak{g}^\ast_{\text{reg}})
\le 3
$. If $\codim_{\frak{g}^\ast}(\frak{g}^\ast-\frak{g}^\ast_{\text{reg}})
= 3
$ then $\frak{g}^\ast-\frak{g}^\ast_{\text{reg}}$ is purely
of codimension three and is precisely equal to the non-smooth locus of
$\pi$. 
\end{proposition}
\begin{proof} Since $X=X_{\sm}$ \eqref{ref-5.6-35} yields a complex of vector
  bundles on $Y$
\begin{equation}
\label{ref-5.9-43}
0\r \pi^\ast \Omega_{X}(1)\xrightarrow{d\pi^\ast} 
\Oscr_{Y}\otimes_k \frak{g}
\xrightarrow{\rho} (
\Oscr_{Y}\otimes\frak{g}^\ast
)(1) \xrightarrow{d\pi}
\pi^\ast T_{X} \r 0
\end{equation}
which is exact in $\pi^\ast T_{X}$ at the smooth locus of $\pi$
(denoted by $U$ above). Furthermore the locus in $U$ where it is exact in
$(\Oscr_{Y}\otimes\frak{g}^\ast )(1)$ is $W\overset{\text{Prop.\
    \ref{ref-5.2-31}}}{=}Y-Y'=\frak{g}^\ast_{\text{reg}}$. 

Assume $ \codim_{\frak{g}^\ast} (\frak{g}^\ast-
\frak{g}^\ast_{\text{reg}})\ge 4 $.  
Hence the locus where
\eqref{ref-5.9-43} is not exact has codimension $\ge 4$. 
An easy depth computation yields that \eqref{ref-5.9-43} is
exact. Thus $\frak{g}^\ast= \frak{g}^\ast_{\text{reg}}$ and hence $0\in \frak{g}^\ast$ is regular. This is only possible if $\frak{g}$ is
abelian, which we had excluded.

Now assume $ \codim_{\frak{g}^\ast} (\frak{g}^\ast-
\frak{g}^\ast_{\text{reg}})=3 $. Then by a similar depth computation
one finds that \eqref{ref-5.9-43} is exact, except in $\pi^\ast
T_{X}$.  Thus in particular $W=U$, or in other words
$\frak{g}^\ast-\frak{g}^\ast_{\text{reg}}$ coincides with the non-smooth
locus of $\pi$. We find that $\coker d\pi$ is a
Cohen-Macaulay module supported on
$\frak{g}^\ast-\frak{g}^\ast_{\text{reg}}$. Since Cohen-Macaulay
modules have pure support we are done.
\end{proof}

\def\cprime{$'$} \def\cprime{$'$} \def\cprime{$'$}
\ifx\undefined\bysame
\newcommand{\bysame}{\leavevmode\hbox to3em{\hrulefill}\,}
\fi

\end{document}